\documentclass[12pt]{article}

\usepackage{latexsym,amsmath,amsfonts,amssymb}
\usepackage[all]{xy}

\def\square{\hbox{\vrule\vbox{\hrule\phantom{o}\hrule}\vrule}}
\def\endofproofsymbol{\square}

\newenvironment{proof}{\par\noindent\textbf{Proof.}\hskip\labelsep}
                      {\unskip\ \nobreak\hbox{}\hfill
                       \endofproofsymbol{\parfillskip=0pt\par}
                       \vspace{\baselineskip}}

\makeatletter \def\@cite#1#2{{#1\if@tempswa , #2\fi}}
              \def\@biblabel#1{#1.}
\makeatother

\newtheorem{Th}{Theorem}[section]
\newtheorem{Lemma}{Lemma}[section]
\newtheorem{Remark}{Remark}[section]
\newtheorem{df}{Definition}[section]
\newtheorem{Example}{Example}[section]
\newtheorem{Proposition}{Proposition}[section]

\newcommand{\R}{\mathbb{R}}

\begin{document}

\title{A Computational Approach to Essential\\ and Nonessential Objective
Functions\\ in Linear Multicriteria Optimization$^1$}

\author{\sc{A. B. Malinowska}$^2$
\sc{and} \sc{D. F. M. Torres}$^3$}

\date{}

\maketitle


\footnotetext[1]{To be partially presented at
the 23rd IFIP TC 7 International Conference
on System Modelling and Optimization, Cracow,
Poland, July 23-27, 2007. Work supported by
KBN under Bialystok Technical University grant W/WI/17/07;
and the R\&D unit CEOC of the University of Aveiro through
FCT and FEDER/POCI 2010.}

\footnotetext[2]{Assistant Professor, Faculty of Computer Science,
Technical University of Bia\l ystok, Bia\l ystok, Poland}

\footnotetext[3]{Associate Professor, Department of Mathematics,
University of Aveiro, Aveiro, Portugal.}



\paragraph*{Abstract.}
The question of obtaining well-defined criteria for multiple
criteria decision making problems is well-known. One of the
approaches dealing with this question is the concept of nonessential
objective function. A certain objective function is called
nonessential if the set of efficient solutions is the same both with
or without that objective function. In this paper we put together
two methods for determining nonessential objective functions. A
computational implementation is done using a computer algebra
system.

\paragraph*{Key~Words.} Multiobjective optimization problems.
Efficient (Pareto optimal) solutions. Essential/nonessential
objective functions.




\section{Introduction}

Multiple criteria decision making (MCDM) arise in connection with
the solution of problems in the areas of economy, environment,
business and engineering (see \textrm{e.g.}
Refs.~\cite{Zbl0462.90054}-\cite{Zbl0588.90019}). 
Multiobjective programming is concerned with the generation of solution
sets for multiobjective problems that usually include a large or an
infinite number of points, known as efficient solutions. Those
efficient points are then the candidates for an optimal solution of the
MCDM problem.

The question of obtaining well-defined criteria for MCDM problems
is well-known. Often mathematical models done by inexperienced
practitioners lead to redundant formulations, which are not only
deceptive but also computationally cumbersome
(see Refs.~\cite{b5}-\cite{b7} and therein). 
Sometimes the decision-maker end up without any decision support,
while a simple reformulation of the problem would achieve the
desired result. One of the approaches dealing with the question of
obtaining well-defined criteria for MCDM is based on the concept
of nonessential objective function. A certain objective function
is called nonessential if the set of efficient solutions is the
same both with or without that objective function. Information
about nonessential objectives helps a decision maker to get
insights and better understand a problem, and this might be a good
starting point for further investigations or revision of the
model. Dropping out nonessential functions leads to a problem with
a smaller number of objectives, which can be solved more easily.
For this reason, the identification of nonessential objectives is
an important feature in the analysis of multiple criteria
programs.

The seminal papers by Gal and Leberling (Refs.~\cite{b3,b4})
define and investigate nonessential objectives in linear
multiobjective optimization problems; Gal and Hanne
(Refs.~\cite{b6,b7}) study the consequences of dropping
nonessential objectives functions in the application of MCDM
methods. Recently, the concept of nonessential objective has been
generalized by Malinowska to convex multiobjective optimization
problems, and new definitions of weakly nonessential and properly
nonessential objective functions were introduced and investigated
(Refs.~\cite{b10,b11}); a new method to determine if a given
objective function of a certain linear problem is essential or not
has been proved (Ref.~\cite{b12}).

Here we put together, in a constructive and algorithmic way, the
two available methods (Refs.~\cite{b3,b12}) for determining
nonessential objective functions. A computational implementation
is done using the computer algebra system Maple. The plan of the
paper is as follows. Section~\ref{sec:PrelNot} gives the necessary
preliminaries and provides the notation used in the text. In
Section~\ref{sec:3} we develop the theory of nonessential
objectives. Main result of the paper appears in
Section~\ref{sec:4}: the algorithm to determine if a given
objective function of a linear multiobjective problem is essential
or not. Finally, in Section~\ref{sec:examples} we provide some
examples that show the applicability of our methodology and the
convenience of the developed computer software. The paper ends
with some conclusions, the references, and an appendix with all
the Maple code that implements the proposed algorithm.


\section{Preliminaries and Notation}
\label{sec:PrelNot}

In this section we present some general concepts and notations. We
use superscripts for vectors (for example $x^{1}$, or simply $x$
when no confusion can arise), and subscripts for components of
vectors (for example $x_{1}$). All the vectors are assumed to be
column vectors. The symbol $\textbf{1}$ stands for the vector
$[1,\ldots,1]^{T}$. For two vectors $x,x'\in \R^{k}$ we define the
relations (Ref.~\cite{b14}):
\begin{equation*}
    x\geqq x'\Leftrightarrow \forall i\in \{1,\ldots,k\}: x_{i}\geq
    x'_{i} \, ,
\end{equation*}
\begin{equation*}
    x\geq x'\Leftrightarrow \forall i\in \{1,\ldots,k\}: x_{i}\geq x'_{i}
    \wedge \exists i\in \{1,\ldots,k\}: x_{i}> x'_{i} \, ,
\end{equation*}
\begin{equation*}
    x> x'\Leftrightarrow \forall i\in \{1,\ldots,k\}: x_{i}>
    x'_{i} \, .
\end{equation*}

Throughout this paper we study the following linear multiobjective
optimization problem:
\begin{equation}\label{p}
\max \{F^{n+1}(x): x\in X\},
\end{equation}
where
\begin{equation*}
 X=\{x\in \R^{k}: Ax\leqq b\, , x\geqq 0\}\, , \quad
 A\in \R^{m\times k}\, , \quad b\in \R^{m}
\end{equation*}
is the feasible set, and
\begin{equation*}
 F^{n+1}(x)=Cx=[(c^{1})^{T}x,\ldots,(c^{n+1})^{T}x]^{T}\, ,
\quad c^{i}\in\R^{k} (i=1,\ldots,n+1)\, , \quad n\geq 1
\end{equation*}
is the vector of objective functions: $f_{i}:\R^{k}\rightarrow \R$
$(i=1,\ldots,n+1)$. We are using ``$\max$'' to mean that we want to
maximize all the objective functions simultaneously. This involves
no loss of generality. In general it does not exist a solution that
is optimal with respect to every objective function, and one is lead
to the concept of Pareto optimality.
\begin{df}
A vector $x^{0}\in X$ is said to be an efficient (Pareto optimal)
solution of problem (\ref{p}) if there exists no $x\in X$ such that
$F^{n+1}(x)\geq F^{n+1}(x^{0})$. The set of efficient solutions of
problem (\ref{p}) is denoted by $X^{n+1}_{E}$.
\end{df}


\section{Nonessential Objectives}
\label{sec:3}

Let $X^{n}_{E}$ denote the set of efficient solutions of problem
(\ref{p}) without one objective function $f_{k}$, $k\in
\{1,\ldots,n+1\}$. Without loss of generality we assume $k=n+1$.
\begin{df}
The objective function $f_{n+1}$ is said to be nonessential in
problem (\ref{p}) if $X^{n}_{E}=X^{n+1}_{E}.$ An objective function
which is not nonessential is called essential.
\end{df}
We now recall two theorems that characterize a nonessential
objective function, and which will be used in the proposed
algorithm to determine if a given objective function is essential
or not.
\begin{Th} (Ref.~\cite{b3})
\label{GL} The objective function $f_{n+1}$ is nonessential in
(\ref{p}) if the following holds:
\begin{equation}\label{wlg}
    c^{n+1}=\sum_{i=1}^{n}\alpha_{i}c^{i} \, , \quad \alpha_{i}\geq
    0 \quad (i=1,\ldots,n) \, .
\end{equation}
\end{Th}
Theorem~\ref{GL} is very useful because condition (\ref{wlg}) is
easy to check. Unfortunately, it is only a sufficient condition.
Theorem~\ref{AM} below gives a necessary and sufficient condition
for an objective function to be nonessential.

Let $X_{n+1}$ denote the set of solutions of the single objective
optimization problem
\begin{equation}\label{}
   \max \{f_{n+1}(x): x\in X\} \, .
\end{equation}
In other words,
\begin{equation}\label{1}
    X_{n+1}=\{x^{0} \in X:\forall x\in X\
  f_{n+1}(x^{0})\geq f_{n+1}(x)\}.
\end{equation}

\begin{Th} (Refs.~\cite{b10,b11})
\label{AM} If the set $X$ is nonempty and bounded (in other words,
if $X$ is a convex polyhedron), then the objective function $f_{n+1}$
is nonessential in (\ref{p}) if and only if the following three
conditions hold:
\begin{itemize}
\item[\text{(i)}] $\forall x \in X\backslash X^{n}_{E}$
$\exists x'\in \R^{k} : F^{n+1}(x')\geq F^{n+1}(x)$;
\item[\text{(ii)}] $X^{n}_{E}\cap X_{n+1}\neq\emptyset$;
\item[\text{(iii)}] $X^{n}_{E}\subset X^{n+1}_{E}$.
\end{itemize}
\end{Th}


\section{Main Results}
\label{sec:4}

The theory described in the previous section enable us to work out
on a computational algorithm to test if an objective function
$f_{n+1}$ in problem (\ref{p}) is essential or not. The proposed
algorithm consists of eight steps and each one has been
implemented in the computer algebra system Maple (see Appendix on
page~\pageref{app}).

\medskip

\noindent \textbf{Step~0.} Does
$c^{n+1}=\sum_{i=1}^{n}\alpha_{i}c^{i}$,
$\alpha_{i}\geq 0$ $(i=1,\ldots,n)?$\\
If the answer is "TRUE", then the objective function $f_{n+1}$ is
nonessential by Theorem~\ref{GL}. Otherwise, we go to Step~1.

\medskip

We implement Step~0 as a Maple command \texttt{glm} (Gal-Leberling
method), which receives in its first argument a list with the
objective functions, and in its second argument the number of
variables of the problem.

\begin{Example}\label{e0.1}
Consider the multiobjective optimization problem
\begin{equation*} \max \{F^{4}(x): x\in X\},
\end{equation*}
where
\begin{equation*}
\begin{split}
  f_{1}(x)&=x_{1}+3x_{2}, \\
  f_{2}(x)&=3x_{1},\\
  f_{3}(x)&=2x_{1}+x_{2},\\
  f_{4}(x)&=-3x_{1}-x_{2},
\end{split}
\end{equation*}
and $X\subset \R^{2}$. With our Maple package we do:
\begin{verbatim}
> glm([x1+3*x2,3*x1,2*x1+x2,-3*x1-x2],2);
            false
\end{verbatim}
We conclude that the answer to Step~0 is "FALSE", so
we go to Step~1.\\
Now, let us change the order of objective functions as follows:
\begin{equation*}
\begin{split}
  f_{1}(x)&=x_{1}+3x_{2}, \\
  f_{2}(x)&=3x_{1},\\
f_{3}(x)&=-3x_{1}-x_{2},\\
 f_{4}(x)&=2x_{1}+x_{2}.
\end{split}
\end{equation*}
This time our \texttt{glm} procedure gives
\begin{verbatim}
> glm([x1+3*x2,3*x1,-3*x1-x2,2*x1+x2],2);
            true
\end{verbatim}
meaning that the answer to Step~0 is "TRUE". Thus, one can
conclude that the objective function $f_{4}(x)=2x_{1}+x_{2}$ is
nonessential.
\end{Example}

\begin{Example}\label{e0.2}
Consider the problem:
\begin{equation*} \max \{F^{3}(x): x\in X\},
\end{equation*}
where
\begin{equation*}
\begin{split}
f_{1}(x)&=x_{1}+x_{2}+x_{3}, \\
  f_{2}(x)&=-x_{1}+x_{2}+x_{3},\\
  f_{3}(x)&=x_{1}+x_{2},
\end{split}
\end{equation*}
and $X\subset \R^{3}$. We have:
\begin{verbatim}
> glm([x1+x2+x3,-x1+x2+x3,x1+x2],3);
            false
\end{verbatim}
Answer to Step~0 is "FALSE", and we go to Step~1.
\end{Example}

\medskip

\textbf{Step~1.} We test condition (i) of Theorem~\ref{AM}.

\medskip

Our method is based on the following observations. Let
\begin{equation*}
U=\{x\in \R^{k}:Cx\geq 0\}.
\end{equation*}

\begin{Remark} (Ref.~\cite{b12})
If $U\neq\emptyset$, then condition (i) of Theorem~\ref{AM} holds.
\end{Remark}

\begin{Th} (Ref.~\cite{b8})
\label{TM} A sufficient condition for $X_{E}^{n+1}=X$ is
$U=\emptyset$. If $intX\neq\emptyset$ (where int stands for the
interior of a set), then this condition is also necessary.
\end{Th}
In Step~1 we solve the problem: does $U\neq\emptyset$? If the
answer is "FALSE", then $X_{E}^{n+1}=X$ and we go to Step~2.
Otherwise, we know that condition (i) of Theorem~\ref{AM} holds
and we go to Step~5. In order to verify equality of sets
$U=\emptyset$ we solve the problem:
\begin{equation}\label{s1}
\max \left\{\sum_{i=1}^{n+1}v_{i}:(x,v)\in V\right\},
\end{equation}
where
\begin{equation*}
V=\{(x,v)\in \R^{k+n+1} : -Cx+v=0\, , v\geqq 0\} \, .
\end{equation*}
\begin{Remark} (Ref.~\cite{b8})
One has $U=\emptyset$ if and only if problem (\ref{s1}) has zero
as the optimal objective function value.
\end{Remark}

\begin{Example}\label{e1.1}
Consider the problem from Example~\ref{e0.1}:
\begin{equation*}
\max\{F^{4}(x): x\in X\},
\end{equation*}
where
\begin{equation*}
\begin{split}
  f_{1}(x)&=x_{1}+3x_{2}, \\
  f_{2}(x)&=3x_{1},\\
  f_{3}(x)&=2x_{1}+x_{2},\\
  f_{4}(x)&=-3x_{1}-x_{2},
\end{split}
\end{equation*}
and $X\subset \R^{2}$. Using our Maple command \texttt{step1} we
obtain:
\begin{verbatim}
> step1(2,[x1+3*x2, 3*x1, 2*x1+x2, -3*x1-x2]);
            false
\end{verbatim}
We conclude that Step~1 has answer "FALSE". Thus, $X_{E}^{4}=X$
and we go to Step~2.
\end{Example}

\begin{Example}\label{e1.2}
Consider again the problem from Example~\ref{e0.2}:
\begin{equation*}
\max\{F^{3}(x): x\in X\},
\end{equation*}
where
\begin{equation*}
\begin{split}
f_{1}(x)&=x_{1}+x_{2}+x_{3}, \\
  f_{2}(x)&=-x_{1}+x_{2}+x_{3},\\
  f_{3}(x)&=x_{1}+x_{2},\\
\end{split}
\end{equation*}
and $X\subset \R^{3}$. We obtain
\begin{verbatim}
> step1(3,[x1+x2+x3, -x1+x2+x3, x1+x2]);
            true
\end{verbatim}
so the answer to Step~1 is "TRUE". Thus, condition (i) of
Theorem~\ref{AM} holds and we go to Step~5.
\end{Example}

\medskip

\textbf{Step~2.} Let  $\tilde{U}=\{x\in \R^{k}:\tilde{C}x\geq
0\}$, where $\tilde{C}=[(c^{1})^{T},\ldots,(c^{n})^{T}]^{T}$. In
Step~2 we address the following question: does
$\tilde{U}\neq\emptyset$?

\medskip

The method we use is the same as the one described in Step~1. If
the answer is "FALSE", then Theorem~\ref{TM} implies $X_{E}^{n}=X$
and the objective function $f_{n+1}$ is nonessential. Otherwise,
we go to Step~3.

\begin{Example}\label{e2}
Let us consider again the problem from Example~\ref{e1.1}. Using
our Maple command \texttt{step2} one has:
\begin{verbatim}
> step2(2,[x1+3*x2, 3*x1, 2*x1+x2, -3*x1-x2]);
            true
\end{verbatim}
Since the answer to Step~2 is "TRUE" we go to Step~3.
\end{Example}

\medskip

\textbf{Step~3.} In this step we solve the problem: does
$intX\neq\emptyset$?

\medskip

Our method is based on the following remark.
\begin{Remark}
If $intX\neq\emptyset$, then problem $\max\{a:(x,v,a)\in V\}$,
where
\begin{equation*}
V=\left\{(x,v,a)\in \R^{k+m+1} : Ax+v+a\textbf{1}=b\, , v\geqq
0,x\geqq \varepsilon\textbf{1}\, , a\geq 0\right\}
\end{equation*}
with $\varepsilon >0$ sufficiently small, has an optimal objective
function value greater than zero.
\end{Remark}
If in Step~3 the answer is "TRUE", then Theorem~\ref{TM} implies
$X_{E}^{n}\neq X$, and the objective function $f_{n+1}$ is
essential. Otherwise, we go to Step~4.

In order to use the simplex package already available from the
Maple system, we put $\varepsilon =0,001$. We note that by default
we are assuming that all $x$ variables are greater or equal than
zero (the user does not need to mention this explicitly in the
definition of the set $X$ while using our Maple package).

\begin{Example}\label{e3}
Let us continue the problem from Examples~\ref{e1.1} and \ref{e2}
with
\begin{equation*}
    X=\{x\in \R^{2}:x_{1}\leq 1,x_{2}\leq 1,x\geqq 0\}.
\end{equation*}
Using our Maple command \texttt{step3} we obtain
\begin{verbatim}
> step3(3,{x1 <= 1, x2 <= 1});
            true
\end{verbatim}
Since the answer to Step~3 is "TRUE", the objective function
$f_{4}(x)=-3x_{1}-x_{2}$ is essential.\\
Now, let us consider a different problem by changing the set $X$
as follows:
\begin{equation*}
    X=\{x\in \R^{2}:x_{1}+x_{2}\leq 1,-x_{1}-x_{2}\leq -1,x\geqq 0\}.
\end{equation*}
Now we obtain
\begin{verbatim}
> step3(2,{x1+x2<=1,-x1-x2<=-1});
            false
\end{verbatim}
Since Step~3 has the answer "FALSE", we go to Step~4.
\end{Example}

\medskip

\textbf{Step~4.} In this step we solve the problem: does
$X_{E}^{n}=X$?

\medskip

Our method is the following: we compute a vertex $x^{0}$ of $X$
and test if $x^{0}$ is an element of $X_{E}^{n}$.

\begin{Th}\label{B} (Ref.~\cite{b1})
Let $x^{0}\in X$ be given. Solve the problem
\begin{equation}\label{s4}
\max\left\{\sum_{i=1}^{n}\epsilon_{i}:(x,\epsilon)\in S\right\}
\end{equation}
with
\begin{equation*}
S=\{(x,\epsilon)\in \R^{k+n}:x\in X,f_{i}(x)-\epsilon
_{i}=f_{i}(x^{0}),\epsilon _{i}\geq 0,i=1,\ldots,n\} \, .
\end{equation*}
The vector $x^{0}$ is efficient if and only if problem (\ref{s4})
has zero as the optimal objective function value.
\end{Th}

If $x^{0}$ is not efficient, then the answer from Step~4 is
"FALSE" and the objective function $f_{n+1}$ is essential.
Otherwise, we compute the next vertex of $X$ and check if it is
efficient or not. Our procedure stops as soon as a non-efficient
vertex is found.

\begin{Example}\label{e4.2}
Consider the problem:
\begin{equation*}
\max\{F^{3}(x): x\in X\},
\end{equation*}
where
\begin{equation*}
\begin{split}
f_{1}(x)&=x_{1}+x_{2} ,\\
f_{2}(x)&=x_{1},\\
f_{3}(x)&=-3x_{1}-x_{2},
\end{split}
\end{equation*}
and
\begin{equation*}
X=\{x\in \R^{2}: x_{1}+x_{2}\leq 1,-x_{1}-x_{2}\leq -1,x\geqq 0\}.
\end{equation*}
Answers from Steps~1 to 4 are easily obtained from the respective
commands of our Maple package:
\begin{verbatim}
> step1(2, [x1 + x2, x1, -3 x1 - x2]);
            false
> step2(2, [x1 + x2, x1, -3 x1 - x2]);
            true
> step3(2,{x1+x2<=1,-x1-x2<=-1});
            false
> step4([x1+x2,x1,-3*x1-x2],{x1+x2<=1,-x1-x2<=-1},2);
            false
\end{verbatim}
We conclude that $f_{3}(x)=-3x_{1}-x_{2}$ is essential and that
$X_{E}^{2}\subset X_{E}^{3}$.
\end{Example}

In the case all vertices are efficient, two situations may appear:
objective function $f_{n+1}$ may be essential (answer "FALSE") or
not (answer "TRUE"). To distinguish between these two cases  we
apply the following remark.

\begin{Remark}\label{s4.1}
Let $x^{1},x^{2},\ldots,x^{p}$ be all vertexes of $X$ and
$intX=\emptyset$. Then, $X_{E}^{n}=X$ if and only if there exists
a vector $w>0$ with $\sum _{i=1}^{n}w_{i}=1$ such that
$w^{T}F^{n}(x^{1})=w^{T}F^{n}(x^{2})=\cdots=w^{T}F^{n}(x^{p})$.
\end{Remark}
\begin{proof}
As far as $intX=\emptyset$, we have $X\subset H$, where $H$ is a
hyperplane. Therefore, $X_{E}^{n}=X$ if and only if there exists
$w>0$ $(\sum _{i=1}^{n}w_{i}=1)$ such that for all $x$ in $X$, $x$
is a solution of the problem $\max\{w^{T}F^{n}(x): x\in X\}$ (see
for instance Ref.~\cite{b2}, p.~54). This completes the proof.
\end{proof}

In order to use the simplex method as provided by Maple, we change
condition from Remark~\ref{s4.1} into the form
\begin{equation*}
\exists w\in \R^{n} :
w^{T}F^{n}(x^{1})=w^{T}F^{n}(x^{2})=\cdots=w^{T}F^{n}(x^{p}), \sum
_{i=1}^{n}w_{i}=1,w^{T}\textbf{1}\geq \varepsilon \textbf{1},
\end{equation*}
where $\varepsilon>0$ is sufficiently small. In our procedure we set
$\varepsilon=0,00001$.

\begin{Example}\label{e4.1}
Let us continue the problem from Examples~\ref{e1.1}, \ref{e2} and
\ref{e3} with
\begin{equation*}
    X=\{x\in \R^{2}:x_{1}+x_{2}\leq 1,-x_{1}-x_{2}\leq -1,x\geqq 0\}.
\end{equation*}
Our Maple command \texttt{step4} give us
\begin{verbatim}
> step4([x1+3*x2,2*x1+x2,3*x1,-3*x1-x2],{x1+x2<=1,-x1-x2<=-1},2);
  true
\end{verbatim}
Thus, $X_E^2 = X$ and it follows that $f_{4}(x)=-3x_{1}-x_{2}$ is
nonessential.
\end{Example}

Now we show an example where all vertices are efficient but $X_E^n
\ne X$.

\begin{Example}
Let us consider the problem
\begin{equation*}
\max\{F^{3}(x): x\in X\},
\end{equation*}
where
\begin{equation*}
\begin{split}
f_1(x_1,x_2,x_3) &= - x_1 - 2 x_2 + 2 x_3, \\
f_2(x_1,x_2,x_3) &= 2 x_1 + 3 x_2, \\
f_3(x_1,x_2,x_3) &= - x_1 - x_2 - 2 x_3,
\end{split}
\end{equation*}
and
\begin{multline*}
X=\Bigl\{x\in \R^{3}: x_2+x_3 \le 2,-x_2-x_3 \le -2,x_1+x_2+x_3 \le 3, \\
-x_1-x_2-x_3 \le -2,x_1+x_2 \le 2 \, , x \geqq 0 \Bigr\} \, .
\end{multline*}
Using our Maple commands we obtain:
\begin{verbatim}
> step1([-x1-2*x2+2*x3, 2*x1+3*x2,-x1-x2-2*x3]);
  false

> step2([-x1-2*x2+2*x3, 2*x1+3*x2,-x1-x2-2*x3]);
  true

> step3(3,{x2+x3<=2,-x2-x3<=-2,x1+x2+x3<=3,-x1-x2-x3<=-2,x1+x2<=2});
  false

> step4([-x1-2*x2+2*x3,2*x1+3*x2,-x1-x2-2*x3],
        {x2+x3<=2,-x2-x3<=-2,x1+x2+x3<=3,-x1-x2-x3<=-2,x1+x2<=2},3);
  false
\end{verbatim}
We conclude that the objective function $f_3$ is essential and
that $X_E^2 \subset X_E^3$.
\end{Example}

\medskip

\textbf{Step~5.} In this step we calculate all vertexes of
$X_{n+1}$ (see (\ref{1})).

\medskip

Let $X^{W}_{n+1}=\{x^{1},x^{2},\ldots,x^{q}\}$ be the set of all
vertexes of $X_{n+1}$. We have:
\begin{equation*}
X_{n+1}=\left\{x\in \R^{k}: x=\sum _{j=1}^{q}\alpha _{j}x^{j}\, ,
\sum _{j=1}^{q}\alpha _{j}=1\, ,
\alpha_{j}\geq0,j=1,\ldots,q\right\} \, .
\end{equation*}

\begin{Example}\label{e5}
Consider again the problem from Example \ref{e1.2} with
\begin{equation*}
    X=\{x\in \R^{3}:x_{1}\leq 1,x_{2}\leq 1,x_{3}\leq 1,x\geqq 0\}.
\end{equation*}
Our procedure \texttt{step5} give us
\begin{verbatim}
> step5([x1+x2+x3,-x1+x2+x3,x1+x2],{x1<=1,x2<=1,x3<=1});
\end{verbatim}
\vspace*{-0.3cm}
$$
\left\{ \left\{ {\it x_1}=1,{\it x_2}=1, {\it x_3}=1 \right\} ,
\left\{ {\it x_1}=1,{\it x_2}=1,{\it x_3}=0 \right\}  \right\}
$$
Hence, $X^{W}_{3}=\{[1,1,0]^{T},[1,1,1]^{T}\}$ and
$$
X_{3}=\{x:x=\alpha [1,1,0]^{T}+(1-\alpha)[1,1,1]^{T},0\leq \alpha
\leq 1\}.
$$
\end{Example}

\medskip

\textbf{Step~6.} In this step we solve the following problem
(condition (ii) of Theorem \ref{AM}): does $X_{n+1}\cap
X^{n}_{E}\neq \emptyset$?

\medskip

The basic idea of our method is to use Theorem~\ref{Y}.
\begin{Th}\label{Y} (Ref.~\cite{b15})
Let $Z=\{x^{1},x^{2},\ldots,x^{q}\}\subset X$. If $Z \subset
X\setminus X^{n}_{E}$, then
\begin{equation*}
\left\{x\in \R^{k}:x=\sum _{j=1}^{q}\alpha _{j}x^{j},\sum
_{j=1}^{q}\alpha _{j}=1,
\alpha_{j}\geq0,j=1,\ldots,q\right\}\subset X\setminus X^{n}_{E}
\, .
\end{equation*}
\end{Th}
Having in mind Theorem~\ref{Y}, it is sufficient to consider only
vertexes of $X_{n+1}$. Applying Theorem~\ref{B} we check if there
exists a vertex of $X_{n+1}$ which belongs to $X^{n}_{E}$. If the
answer is "TRUE", then the condition (ii) of Theorem~\ref{AM}
holds and we go to Step~7. Otherwise, we conclude that the
objective function $f_{n+1}$ is essential.

\begin{Example}\label{e6.1}
Let us  continue the problem from Examples \ref{e1.2} and
\ref{e5}. Our procedure \texttt{step6} give the answer "TRUE"
\begin{verbatim}
> step6([x1+x2+x3,-x1+x2+x3,x1+x2],{x1<=1,x2<=1,x3<=1});
            true
\end{verbatim}
and we proceed to Step~7.
\end{Example}

\begin{Example}\label{e6.2}
Now we consider a problem borrowed from Ref.~\cite{Zbl0588.90019}:
\begin{equation*}
\max\{F^{3}(x): x\in X\},
\end{equation*}
where
\begin{equation*}
\begin{split}
f_{1}(x)&= x_1+2 x_2-x_3+3 x_4+ 2x_5+x_7,\\
f_{2}(x)&= x_2+x_3+ 2 x_4 + 3 x_5+x_6,\\
f_{3}(x)&= x_1 + x_3-x_4-x_6-x_7,
\end{split}
\end{equation*}
and
\begin{multline*}
X=\{x\in \R^{7}: x_1+2 x_2+x_3+x_4+2 x_5+x_6+2 x_7 \le 16,\\
-2 x_1-x_2+x_4+2 x_5+x_7 \le 16,\\
-x_1+x_3+2 x_5-2 x_7 \le 16, x_2+2 x_3-x_4+x_5-2 x_6-x_7 \le 1,
x\geqq 0\}.
\end{multline*}
With our Maple command \texttt{step6} we obtain
\begin{verbatim}
> step6([x1+2*x2-x3+3*x4+2*x5+x7,x2+x3+2*x4+3*x5+x6,x1+x3-x4-x6-x7],
        {x1+2*x2+x3+x4+2*x5+x6+2*x7<=16,-2*x1-x2+x4+2*x5+x7<=16,
         -x1+x3+2*x5-2*x7<=16,x2+2*x3-x4+x5-2*x6-x7<=1});
  false
\end{verbatim}
and since the answer is "FALSE", we conclude that the objective
function $f_{3}(x) = x_1 + x_3-x_4-x_6-x_7$ is essential.
\end{Example}

\medskip

\textbf{Step~7.} In this step we solve the following problem
(condition (iii) of Theorem~\ref{AM}): does $X_{E}^{n}\subset
X_{E}^{n+1}$?

\medskip

Our method is based on the following observations.
\begin{Proposition} (Refs.~\cite{b10,b9})
\label{Prop:4dot12} If the vector-valued function $F^{n}$ is
one-to-one on the set $X^{n}_{E}$, then condition (iii) of
Theorem~\ref{AM} holds.
\end{Proposition}

Let
\begin{equation*}
    \langle X_{E}^{n}\rangle=\{ x^{i}-x^{j}: x^{i},x^{j}\in X^{n}_{E}\}.
\end{equation*}

\begin{Lemma} (Ref.~\cite{b16})
The vector-valued function $F^{n}$ is one-to-one on the set
$X^{n}_{E}$ if and only if $ KerF^{n}\cap \langle X_{E}^{n}\rangle
=\emptyset$ (Ker stands for the kernel of a map).
\end{Lemma}
In practice it is usually impossible to determine the set $\langle
X_{E}^{n}\rangle$. For this reason, in our Maple procedure we use
the following set:
\begin{equation*}
    \langle X_{WE}^{n}\rangle=\{ x^{i}-x^{0}: x^{i},x^{0}\in X^{n}_{E}\},
\end{equation*}
where $x^{i},x^{0}$ are vertexes of $X^{n}_{E}$ and $x^{0}$ is
free but fixed. Obviously, $\langle X_{E}^{n}\rangle \subset
Lin\{\langle X_{WE}^{n}\rangle \}$.
\begin{Remark}
If $KerF^{n}\cap Lin \{\langle X_{WE}^{n}\rangle\} =\emptyset$,
then condition (iii) of Theorem~\ref{AM} holds (that is,
$X_{E}^{n}\subset X_{E}^{n+1}$).
\end{Remark}

If the answer from Step~7 is "TRUE", then the objective function
$f_{n+1}$ is nonessential. Otherwise, we know that
$X_{E}^{n+1}\subset X_{E}^{n}$.

\begin{Example}
\label{ex:st7}
Let us  continue the problem from Examples \ref{e1.2}, \ref{e5}
and \ref{e6.1}. We obtain:
\begin{verbatim}
> step7([x1+x2+x3,-x1+x2+x3,x1+x2],{x1<=1,x2<=1,x3<=1},3);
            true
\end{verbatim}
The answer from Step~7 is "TRUE", hence the objective function
$f_{3}(x)=x_{1}+x_{2}$ is nonessential.
\end{Example}


\section{Illustrative Examples}
\label{sec:examples}

We have implemented all the steps described in Section~4 together
in a single Maple command called \texttt{nonessential} (see
Figure~\ref{fig}). This main procedure receives a list with the
definition of the objective functions, and the set $X$ of
constraints in the second argument. Below we give some examples of
computer sessions with our Maple package. The interested reader
may download it from
\texttt{[http://www.mat.ua.pt/delfim/essential.html]} and find
there more examples than the ones we are able to provide here. We
invite and welcome the reader to experiment our Maple package with
her/his own problems.

\begin{center}
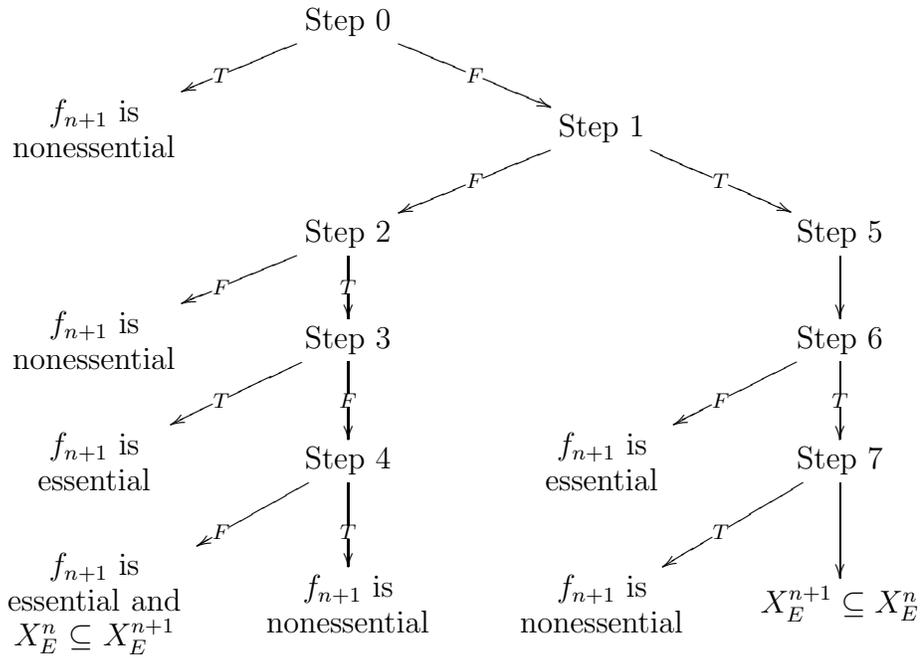
\begin{figure}
\xymatrix @-0.5pc {
  & \ar[dl]|T \text{Step 0}\ar[dr]|F & & \\
\txt<6pc>{$f_{n+1}$ is nonessential} & & \ar[dl]|F\text{Step 1}\ar[dr]|T & \\
 & \ar[dl]|F\text{Step 2}\ar[d]|T & & \text{Step 5}\ar[d] \\
\txt<6pc>{$f_{n+1}$ is nonessential} & \ar[dl]|T\text{Step 3}\ar[d]|F & & \ar[dl]|F\text{Step 6}\ar[d]|T \\
\txt<6pc>{$f_{n+1}$ is essential} & \ar[dl]|F\text{Step 4}\ar[d]|T & \txt<6pc>{$f_{n+1}$ is essential} & \ar[dl]|T\text{Step 7}\ar[d] \\
\txt<6pc>{$f_{n+1}$ is essential and  $X_E^n \subseteq X_E^{n+1}$} & \txt<6pc>{$f_{n+1}$ is nonessential} & \txt<6pc>{$f_{n+1}$ is nonessential} & X_E^{n+1} \subseteq X_E^{n}
}
\caption{Scheme of the method implemented in Maple.}
\label{fig}
\end{figure}
\end{center}

We begin with a simple example with three objective functions and
three variables.

\begin{verbatim}
> nonessential([x1+x2,x1+x2+x3,-3*x1-3*x2-x3],{x1+x2+x3 <= 1});
Objective function -3*x1-3*x2-x3 is essential from Step 3
\end{verbatim}

Next we consider a problem with four objective functions and five
variables. It turns out that the problem may be reduced to a simpler
one with the same set of efficient solutions.

\begin{verbatim}
> st:= time():
  nonessential([x1+x2+x3+x4+x5,-x1+x2+x3+x4+x5,-x1-x2+x3+x4+x5,x1+x2],
               {x1<=1,x2<=1,x3<=1,x4<=1,x5<=1});
  printf("%a seconds\n",time() - st);

Objective function x1+x2 is nonessential from step 7
2.072 seconds
\end{verbatim}

The only situation where our Maple procedure \texttt{nonessential}
can not conclude that a given function is essential or not, is
when one reaches Step~7 and the answer is not true.

\begin{verbatim}
> st:= time():
  nonessential([-x3-x4,-x5-x6,-x4-x6],
    {x1+3*x2<=24,3*x1+x2<=24,x1+4*x2+x3-x4<=40,
     -x1+4*x2-x3+x4<=-40,4*x1+x2+x5-x6<=40,-4*x1-x2-x5+x6<=-40});
  printf("%a seconds\n",time() - st);

X_E^3 C X_E^2 from step 7
7.251 seconds
\end{verbatim}

Our Maple package is useful to identify redundant objective
functions. We finish with an example where the mathematical model
can be simplified by elimination of two of the objective functions.

\begin{verbatim}
> nonessential([x1+3*x2,2*x1+x2,3*x1,-3*x1-x2],
               {x1+x2<=1,-x1-x2<=-1});
  Objective function -3*x1-x2 is nonessential from step 4
> nonessential([x1+3*x2,2*x1+x2,3*x1],{x1+x2<=1,-x1-x2<=-1});
  Objective function 3*x1 is nonessential from step 7
> nonessential([x1+3*x2,2*x1+x2],{x1+x2<=1,-x1-x2<=-1});
  Objective function 2*x1+x2 is essential from step 6
> nonessential([2*x1+x2,x1+3*x2],{x1+x2<=1,-x1-x2<=-1});
  Objective function x1+3*x2 is essential from step 6
\end{verbatim}


\section{Conclusion}

There are theoretical and practical reasons for developing a
method to find nonessential objective functions. In this paper we
present such a method and its implementation in Maple. Our
algorithm is based on necessary and sufficient conditions for an
objective function to be nonessential, and need only to solve a
finite number of single objective linear optimization problems.
Examples showing the usefulness of our Maple package are provided:
identification of nonessential objective functions permits to
simplify the correspondent mathematical model.



\section*{Appendix -- Maple Definitions}
\label{app}

Our Maple definitions follow below. The reader can download the
code from [\texttt{http://www.mat.ua.pt/delfim/essential.html}]
together with many more examples than the ones we are able to
provide in the paper.

\smallskip

We begin by implementing the Gal-Leberling method (see
Examples~\ref{e0.1} and \ref{e0.2}), which is based on the results
of Ref.~\cite{b3}.

\footnotesize
\begin{verbatim}
> ##################################################################
> # GL method; returns true if F[-1] is nonessential;
> # returns false if one can not conclude nothing from GL method
> ##################################################################
> glm := proc(F,NumVar)
>   local c, f, v, LV, lc, SolSet, SS, i;
>   c := proc(var,exp)
>     local v;
>     v:= select(has,exp+abm,var);
>     if v = NULL then
>       return(0);
>     else
>       return(v/var):
>     fi:
>   end proc;
>   f := o -> if type(o,numeric) then o else 0 fi:
>   v := (exp,n) -> Vector([seq(f(c(x||i,exp)),i=1..n)]):
>   LV := [seq(LinearAlgebra[VectorScalarMultiply](v(F[i],NumVar),
               alpha||i),i=1..nops(F)-1)];
>   lc := Vector(NumVar);
>   for i in LV do
>     lc := LinearAlgebra[VectorAdd](lc,i):
>   od;
>   SolSet := solve({seq(lc[i]=c(x||i,F[-1]),i=1..NumVar)});
>   if SolSet = NULL then
>     return(false);
>   else
>     SS := {seq(simplex[maximize](alpha||i,SolSet,
                      NONNEGATIVE),i=1..nops(F)-1)};
>     if SS = {{}} then
>       return(false);
>     else
>       return(true);
>     fi:
>   fi:
> end proc:
\end{verbatim}
\normalsize

For illustrative examples on how to use the procedures
\texttt{step1} and \texttt{step2} see Examples~\ref{e1.1},
\ref{e1.2} and \ref{e2}.

\footnotesize
\begin{verbatim}
> ###################
> # step 1
> ###################
> step1 := proc(F)
>   local of, const, S;
>   of := add(v||i,i=1..nops(F));
>   const := seq(-F[i]+v||i=0,i=1..nops(F));
>   const := [const,seq(v||i>=0,i=1..nops(F))];
>   S := simplex[maximize](of,const);
>   return(not(evalb(subs(S,of)=0))):
> end proc:
> ###################
> # step 2
> ###################
> step2 := F -> step1(F[1..-2]):
\end{verbatim}
\normalsize
Follows our implementation in Maple for Step~3 (see Example~\ref{e3}).
\footnotesize
\begin{verbatim}
> step3 := proc(NumVar,X)
>   local LHS, RHS, SC1, SC2, SC3, SC, SS3, i, a, epsilon, mylhs, myrhs;
>   epsilon := 0.001;
>   mylhs := E -> if type(lhs(E),numeric) then rhs(E) else lhs(E) fi:
>   myrhs := E -> if type(rhs(E),numeric) then rhs(E) else lhs(E) fi:
>   LHS := [seq(mylhs(i),i=X)];
>   RHS := [seq(myrhs(i),i=X)];
>   SC1 := {seq(LHS[i]+a+v||i = RHS[i],i=1..nops(LHS))};
>   SC2 := {seq(v||i >= 0,i=1..nops(LHS)), a>=0};
>   SC3 := {seq(x||i>=0.001,i=1..NumVar)};
>   SC := SC1 union SC2 union SC3;
>   SS3 := simplex[maximize](a,SC);
>   assign(select(has,SS3,a));
>   if a = 0 then return(false) else return(true) fi;
> end proc:
\end{verbatim}
\normalsize In Step~4 we use a rank method, computing the rank of a
matrix $A$ by the \texttt{Rank} command from the standard
\texttt{LinearAlgebra} package of Maple system. We notice that the
procedure \texttt{step4} does not use the last objective function
(the last objective is given in Maple by \texttt{F[-1]}, and we
exclude it from consideration by writing \texttt{F[1..-2]}). The
auxiliary procedure \texttt{matrixA} is used both by Steps~4 and 7.
The procedure \texttt{Proposition4dot12} is our Maple definition for
Remark~\ref{s4.1}. \footnotesize
\begin{verbatim}
> matrixA := proc(X,NumVar)
>   local c, LHS, row, mylhs;
>   mylhs := E -> if type(lhs(E),numeric) then rhs(E) else lhs(E) fi:
>   c := (var,exp) -> if evalb({select(has,exp,var)} = {}) then
                        0
                      else
                        select(has,exp,var)/var
                      fi:
>   LHS := [seq(mylhs(i)+abm,i=X)];
>   row := (exp,NumVar) -> map(c,[seq(x||i,i=1..NumVar)],exp):
>   return(Matrix(map(row,LHS,NumVar)));
> end proc:
>
> Proposition4dot12 := proc(F,X,LE)
>   local SE, v, ETS, solW,  SS, SV;
>   SE := NULL;
>   for v in LE do
>     SV := subs(v,F);
>     SE := SE, add(SV[i]*w||i,i=1..nops(SV));
>   od;
>   ETS := seq(SE[1]=i,i=SE[2..-1]), add(w||i,i=1..nops(SV))=1;
>   solW := solve({ETS, add(w||i,i=1..nops(SV))=1})
                   union {seq(w||i>=0.00001,i=1..nops(SV))};
>   SS := {seq(simplex[maximize](w||i,solW),i=1..nops(SV))};
>   return(not(remove(i->i={},SS) = {}));
> end proc:
> step4 := proc(F,X,NumVar)
>   local NX,NNumVar,b,A,rnk,dif,LP,zero,efficient,v,admissible,
>         vs,AM,sol,of,cstEps,cst,SC,p,myrhs,mylhs,LE, tv, val;
>   myrhs := E -> if type(rhs(E),numeric) then rhs(E) else lhs(E) fi:
>   mylhs := E -> if type(lhs(E),numeric) then rhs(E) else lhs(E) fi:
>   tv := (x,s) -> myrhs(op(select(has,s,x))):
>   b := X -> Vector([seq(myrhs(i),i=X)]):
>   NX := {seq(mylhs(X[i])+x||(NumVar+i)=myrhs(X[i]),i=1..nops(X))};
>   NNumVar := NumVar+nops(X);
>   A := matrixA(NX,NNumVar);
>   rnk := LinearAlgebra[Rank](A);
>   dif := NNumVar - rnk;
>   LP := combinat[choose]([seq(x||i,i=1..NNumVar)],dif);
>   zero := L -> {seq(i=0,i=L)}:
>   efficient := true;
>   v := Vector([seq(x||i,i=1..NNumVar)]);
>   admissible := sc -> not(member(false,{seq(evalb(myrhs(i)>=0),i=sc)})):
>   LE := NULL;
>   for p in LP while efficient do
>     vs := subs({seq(i=0,i=p)},v);
>     AM := LinearAlgebra[MatrixVectorMultiply](A,vs);
>     sol := solve({seq(AM[i]=b(NX)[i],i=1..nops(NX))});
>     if not(sol = NULL) and admissible(sol) then
>       sol := sol union zero(p):
>       of := add(epsilon||i,i=1..nops(F)-1);
>       cstEps := seq(epsilon||i>=0,i=1..nops(F)-1);
>       cst := seq(F[i]-epsilon||i=subs(sol,F[i]),i=1..nops(F)-1);
>       SC := {cst} union {cstEps} union NX;
>       efficient := evalb(subs(simplex[maximize](of,SC,NONNEGATIVE),of)=0);
>       if efficient then
>         val := [seq(tv(x||i, sol),i=1..NumVar)];
>         if not(member(false,map(i->type(i,numeric),val))) then
>           LE := LE, {seq(x||i=val[i],i=1..NumVar)};
>         fi;
>       fi;
>     fi:
>   od;
>   if not(efficient) then
>     return(efficient);
>   else
>     Proposition4dot12(F[1..-2],X,{LE});
>   fi:
> end proc:
\end{verbatim}
\normalsize
The procedure \texttt{step5} makes use of
an auxiliar procedure \texttt{vert} that
receives three arguments: one solution given
by the simplex method (denoted by \texttt{Sol});
an objective function \texttt{of}; and a set of constraints \texttt{X}.
We remark that in Step~5 only the last objective function is used (that
is given in Maple by \texttt{F[-1]}, where \texttt{F}
is the list of all the objectives under consideration).
\footnotesize
\begin{verbatim}
> vert := proc(Sol,of,X)
>   local v,S,LFV,LS,i,tv,Min,Max,LL,LV,gaa,ga,VFV,freeVar,Sub,LSub,
          s,delFreeVar,varSol,varOF1,varOF,VerifySol,aux,NX, mylhs, myrhs;
>   mylhs := E -> if type(lhs(E),numeric) then rhs(E) else lhs(E) fi:
>   myrhs := E -> if type(rhs(E),numeric) then rhs(E) else lhs(E) fi:
>   v := subs(Sol,of);
>   varOF1 := t -> op(select(i->not(type(i,numeric)),[op(t)])):
>   varOF := of -> map(varOF1,{op(of)}):
>   S := solve(of = v,varOF(of));
>   freeVar := SS -> {seq(mylhs(i),i=select(i->mylhs(i)=myrhs(i),SS))}:
>   varSol := E -> {seq(mylhs(i),i=E)}:
>   LFV := freeVar(S) union (varSol(Sol) minus varOF(of));
>   tv := (x,s) -> myrhs(op(select(has,s,x))):
>   Min := (x,X) -> simplex[minimize](x,X union S,NONNEGATIVE):
>   Max := (x,X) -> simplex[maximize](x,X union S,NONNEGATIVE):
>   LL := [seq([tv(x,Min(x,X)),tv(x,Max(x,X))],x=LFV)];
>   LV := n -> [seq(i[j],j=1..n)]:
>   gaa := (n,m,L) -> if m=n then
                        seq(LV(m),i[m]=L[m])
                      else
                        seq(gaa(n,m+1,L),i[m]=L[m])
                      fi:
>   ga := L -> gaa(nops(L),1,L):
>   VFV := {ga(LL)};
>   Sub := (C1,C2) -> seq({seq(C2[j]=i[j],j=1..nops(i))},i=C1):
>   LSub := Sub(VFV,LFV);
>   delFreeVar := SS -> SS minus {seq(i,i=select(i->mylhs(i)=myrhs(i),SS))}:
>   NX := X union {seq(i>=0,i=varSol(Sol))}:
>   VerifySol := PS -> not(member(false,{seq(evalb(subs(PS,i)),i=NX)})):
>   aux := {seq(subs(s,delFreeVar(S)) union s,s={LSub})}:
>   return(select(VerifySol,aux));
> end proc:
>
> step5 := proc(F,X)
>   local SolSM;
>   SolSM := simplex[maximize](F[-1],X,NONNEGATIVE);
>   return(vert(SolSM,F[-1],X));
> end proc:
\end{verbatim}
\normalsize
Examples~\ref{e6.1} and \ref{e6.2} illustrate the use
of our Maple command \texttt{step6}.
\footnotesize
\begin{verbatim}
> step6 := proc(F,X)
>   local STEP5, of, cstEps, notEfficient, sol, cst, SC;
>   STEP5 := step5(F,X);
>   of := add(epsilon||i,i=1..nops(F)-1);
>   cstEps := seq(epsilon||i>=0,i=1..nops(F)-1);
>   notEfficient := true;
>   for sol in STEP5 while notEfficient do
>     cst := seq(F[i]-epsilon||i=subs(sol,F[i]),i=1..nops(F)-1):
>     SC := {cst} union {cstEps} union X:
>     subs(simplex[maximize](of,SC,NONNEGATIVE),of);
>     notEfficient := evalb(subs(simplex[maximize](of,SC,NONNEGATIVE),of)<>0);
>   end do;
>   return(not(notEfficient));
> end proc:
\end{verbatim}
\normalsize
Follows our Maple definition for Step~7 (see Example~\ref{ex:st7}).
\footnotesize
\begin{verbatim}
> step7a := proc(F,X,NumVar)
>   local b,A,rnk,dif,LP,zero,efficient,v,admissible,
          vs,AM,sol,of,cstEps,cst,SC,p,LE,myrhs;
>   myrhs := E -> if type(rhs(E),numeric) then rhs(E) else lhs(E) fi:
>   b := X -> Vector([seq(myrhs(i),i=X)]):
>   A := matrixA(X,NumVar);
>   rnk := LinearAlgebra[Rank](A);
>   dif := NumVar - rnk;
>   LP := combinat[choose]([seq(x||i,i=1..NumVar)],dif);
>   zero := L -> {seq(i=0,i=L)}:
>   LE := NULL;
>   v := Vector([seq(x||i,i=1..NumVar)]);
>   admissible := sc -> not(member(false,{seq(evalb(myrhs(i)>=0),i=sc)})):
>   for p in LP do
>     vs := subs({seq(i=0,i=p)},v);
>     AM := LinearAlgebra[MatrixVectorMultiply](A,vs);
>     sol := solve({seq(AM[i]=b(X)[i],i=1..nops(X))});
>     if not(sol = NULL) and admissible(sol) then
>       sol := sol union zero(p):
>       of := add(epsilon||i,i=1..nops(F)-1);
>       cstEps := seq(epsilon||i>=0,i=1..nops(F)-1);
>       cst := seq(F[i]-epsilon||i=subs(sol,F[i]),i=1..nops(F)-1);
>       SC := {cst} union {cstEps} union X;
>       efficient := evalb(subs(simplex[maximize](of,SC,NONNEGATIVE),of)=0);
>       if efficient then LE := LE, sol; fi:
>     fi:
>   od;
>   return([LE]);
> end proc:
> ##############################################################
> # In step7b we change X
> # Note: All inequalities must be given in the form Ax <= b
> ##############################################################
> step7b := proc(F,X,NumVar)
>   local TX, SEV, good, sel, mylhs, myrhs;
>   mylhs := E -> if type(lhs(E),numeric) then rhs(E) else lhs(E) fi:
>   myrhs := E -> if type(rhs(E),numeric) then rhs(E) else lhs(E) fi:
>   TX := {seq(mylhs(X[i])+x||(NumVar+i)=myrhs(X[i]),i=1..nops(X))};
>   SEV := step7a(F,TX,NumVar+nops(X));
>   good := (v,nv) -> member(v,{seq(x||i,i=1..nv)}):
>   sel := (es,NumVar) -> select(e->good(mylhs(e),NumVar),es):
>   return(map(sel,SEV,NumVar));
> end proc:
>
> step7 := proc(F,X,NumVar)
>   local C, kernel, S7, tv, five, SD, basis, IBK, myrhs;
>   myrhs := E -> if type(rhs(E),numeric) then rhs(E) else lhs(E) fi:
>   C := matrixA([seq(i=0,i=F[1..-2])],NumVar);
>   kernel := LinearAlgebra[NullSpace](C);
>   if kernel = {} then
>     printf("Objective function %a is nonessential from step 7\n",F[-1]);
>   else
>     S7 := step7b(F,X,NumVar);
>     tv := (x,S) -> myrhs(op(select(has,S,x))):
>     five := (S1,S2,NumVar) -> Vector([seq(tv(x||i,S2)-tv(x||i,S1),i=1..NumVar)]):
>     SD := [seq(five(S7[1],S7[i],NumVar),i=2..nops(S7))];
>     basis := LinearAlgebra[Basis](SD);
>     IBK := LinearAlgebra[IntersectionBasis]([basis,kernel]);
>     if IBK = {} then
>       printf("Objective function %a is nonessential from step 7\n",F[-1]);
>     else
>       printf("X_E^%a C X_E^%a from step 7\n",nops(F),nops(F)-1);
>     fi:
>   fi:
> end proc:
\end{verbatim}
\normalsize
Our Maple procedure \texttt{mm} is based on the theory
introduced in Ref.~\cite{b12} (\texttt{mm} stands for ``Malinowska Method'').
\footnotesize
\begin{verbatim}
> mm := proc(F,X,NumVar)
>   if step1(F) then
>     if step6(F,X) then
>       step7(F,X,NumVar);
>     else
>       printf("Objective function %a is essential from step 6\n",F[-1]);
>     fi;
>   else
>     if not(step2(F)) then
>       printf("Objective function %a is nonessential from step 2\n",F[-1]);
>     else
>       if step3(NumVar,X) then
>         printf("Objective function %a is essential from step 3\n",F[-1]);
>       else
>         if step4(F,X,NumVar) then
>           printf("Objective function %a is nonessential from step 4\n",F[-1]);
>         else
>           printf("Objective function %a is essential from step 4
>                   and X_E^%a C X_E^%a\n",F[-1],nops(F)-1,nops(F));
>         fi:
>       fi:
>     fi:
>   fi:
> end proc:
\end{verbatim}
\normalsize
Follows the main procedure of our Maple package
(see Section~\ref{sec:examples}).
\footnotesize
\begin{verbatim}
> nonessential := proc(F,X)
>   local NumVar, y, cs;
>   cs := x -> [op(x)][-1]:
>   y := sort(map(cs,remove(i->type(i,numeric),
                       map(i->op(i),
                         [seq(op(i),i=F),seq(op(i),i=X)]))))[-1];
>   for NumVar from 1 by 1 while not(evalb(x||NumVar = y)
                                 or evalb(-x||NumVar = y)) do od;
>   if glm(F,NumVar) then
>     printf("Objective function %a is nonessential from GL method\n",F[-1]);
>   else
>     mm(F,X,NumVar);
>   fi:
> end proc:
\end{verbatim}
\normalsize


\end{document}